\documentclass[11pt,amsfonts]{article}
\usepackage{graphicx}
\usepackage{latexsym}
\usepackage{amssymb}
\usepackage{amsmath}
\usepackage{enumerate}

\usepackage{layout}
\usepackage{eufrak}
\newtheorem{prop}{Proposition}

\newtheorem{corollary}{Corollary}
\newtheorem{theorem}{Theorem}
\newtheorem{remark}{Remark}

\def\real{{\mathord{{\rm I\kern-2.8pt R}}}}        
\def\inte{{\mathord{{\rm I\kern-2.8pt N}}}}

\def\sZZ{{\rm Z\kern-2.8ptem{}Z}}

\def\z{{\mathchoice
  {\sZZ}
  {\sZZ}
  {\rm Z\kern-0.30em{}Z}
  {\rm Z\kern-0.25em{}Z} }}
\def\sQQ{{\kern 0.27em \vrule height1.45ex width0.03em depth0em
          \kern-0.30em \rm Q}}
\def\qu{{\mathchoice
    {\sQQ}
    {\sQQ}
  {\kern 0.225em \vrule height1.05ex width0.025em depth0em \kern-0.25em \rm Q}
  {\kern 0.180em \vrule height0.78ex width0.020em depth0em \kern-0.20em \rm Q}
        }}
\def\sCC{{\kern 0.27em \vrule height1.45ex width0.03em depth0em
          \kern-0.30em \rm C}}
\def\complex{{\mathchoice
    {\sCC}
    {\sCC}
  {\kern 0.225em \vrule height1.05ex width0.025em depth0em \kern-0.25em \rm C}
  {\kern 0.180em \vrule height0.78ex width0.020em depth0em \kern-0.20em \rm C}
        }}

\newcommand{\ba}{\begin{array}}
\newcommand{\ea}{\end{array}}
\newcommand{\be}{\begin{equation}}
\newcommand{\ee}{\end{equation}}
\newcommand{\bea}{\begin{eqnarray}}
\newcommand{\eea}{\end{eqnarray}}
\newcommand{\beaa}{\begin{eqnarray*}}
\newcommand{\eeaa}{\end{eqnarray*}}

\newcommand{\eps}{\varepsilon}

\def\z{\zeta}

\font\tenmath=msbm10 \font\sevenmath=msbm7 \font\fivemath=msbm5
\newfam\mathfam \textfont\mathfam=\tenmath
\scriptfont\mathfam=\sevenmath \scriptscriptfont\mathfam=\fivemath

\def \={{\buildrel {\rm (law)} \over =}}

\def\qed{ \hfill \vrule width.25cm height.25cm depth0cm\smallskip}

\newcommand{\basa}{\begin{assumption}}
\newcommand{\easa}{\end{assumption}}

\newcommand{\bas}{\begin{assum}}
\newcommand{\eas}{\end{assum}}

\def\supp{\hbox{\rm supp$\,$}}
\def\span{\hbox{\rm span$\,$}}

\newcommand{\ignore}[1]{}
\author{meryem slaoui }
\date{}
\begin{document}

\renewcommand{\thefootnote}{\fnsymbol{footnote}}

\renewcommand{\thefootnote}{\fnsymbol{footnote}}

\title{Limit behavior of the Rosenblatt Ornstein-Uhlenbeck process with respect to the Hurst index }
\author{Meryem Slaoui and C. A. Tudor\vspace*{0.2in} \\
 $^{1}$ Laboratoire Paul Painlev\'e, Universit\'e de Lille 1\\
 F-59655 Villeneuve d'Ascq, France.\\
\quad meryem.slaoui@math.univ-lille1.fr\\
 \quad tudor@math.univ-lille1.fr\\
\vspace*{0.1in} }

\maketitle

\begin{abstract}
We study the convergence in distribution, as $H\to \frac{1}{2}$ and as $H\to 1$, of the integral $\int_{\mathbb{R}} f(u) dZ^{H}(u) $, where $Z ^{H}$ is a Rosenblatt process with self-similarity index $H\in \left( \frac{1}{2}, 1\right) $ and $f$ is a suitable deterministic function. We focus our analysis on the case of the Rosenblatt Ornstein-Uhlenbeck process, which is the solution of the Langevin equation driven by the Rosenblatt process.
\end{abstract}

\vskip0.3cm

{\bf 2010 AMS Classification Numbers:}   60H05, 60H15, 60G22.

\vskip0.3cm


\section{Introduction}
 The Rosenblatt process is a stochastic process in the second Wiener chaos, i.e. it can be expressed as a multiple integral  of order two with respect to the Wiener process. It is a non-Gaussian self-similar process with stationary increments that exhibits long memory. The Rosenblatt process belongs to the class of so-called Hermite processes, which are self-similar processes with stationary increments in the $q$th Wiener chaos ($q\geq 1$) (the Rosenblatt process is obtained for $q=2$ while for $q=1$ we have the fractional Brownian motion, which is the only Gaussian Hermite process). The Rosenblatt process has been widely studied in the last decades, see e.g. the monographs \cite{PiTa-book} and \cite{T} and the references therein.

There are  several recent research works that investigate the asymptotic  behavior in distribution of some fractional processes (see \cite{BellNu}, \cite{BaiTa}, \cite{ArTu}, \cite{VeTa}) with respect to the Hurst parameter. In particular, in  the case of the Rosenblatt process $(Z ^{H}(t))_{t\geq 0} $ with self-similarity index $H\in (\frac{1}{2}, 1)$, it has been shown in \cite{VeTa} that $ Z^{H}$ converges weakly,   as $H\to \frac{1}{2}$,  in the space of continuous functions $C[0,T]$ (for every $T>0$),   to a Brownian motion while if  $H\to 1$, it tends weakly to the stochastic process  $(t \frac{1}{\sqrt{2}}(Z^{2}-1))_{t\geq 0} $, $Z^{2}-1$ being a so-called {\it centered chi-square random variable}. The case of the generalized Rosenblatt process has been considered in \cite{BaiTa} while the case of the Rosenblatt sheet can be found in \cite{ArTu}. Hermite processes of higher order have been considered in \cite{BellNu}, \cite{ArTu}.

The purpose of this work is to investigate the asymptotic behavior in distribution, with respect to the Hurst parameter, of   the Wiener integral with respect to the Rosenblatt process (or the Wiener-Rosenblatt integral). The Wiener-Rosenblatt  integral  $\int_{\mathbb{R}} f(u) d Z^{H}(u)$ has been introduced in \cite{MaTu}, for a sufficiently regular deterministic function $f$. In a first part, we give the asymptotic behavior in distribution, as $H\to \frac{1}{2}$ and as $H\to 1$,  of the random variable  $\int_{\mathbb{R}} f(u) d Z^{H}(u)$, by assuming suitable integrability condition on $f$. We will then focus on the asymptotic behavior with respect to $H$ of the Rosenblatt Ornstein-Uhlenbeck process (ROU for short) which constitutes the unique solution of the Langevin equation driven by the Rosenblatt process. The ROU process can be expressed in the form of a Wiener-Rosenblatt  integral with a particular kernel $f$. In order to check that this kernel verifies the integrability conditions needed  to apply our general result, we need to use some technical results, in particular the so-called power counting theorem from \cite{TeTa}. We will treat separately the cases of the non-stationary ROU (whose initial value does not depend on $H$) and of the stationary ROU (with initial value depending on the Hurst parameter).  We will prove that the (stationary) ROU converges weakly, as $H\to \frac{1}{2}$, to the (stationary) Gaussian Ornstein-Uhlenbeck process (solution of the Langevin equation driven by the Brownian motion) while as $H\to 1$, the ROU process converges weakly  to  a chi-square random variable  multiplied by a deterministic process. 
Since we deal with processes in the second Wiener chaos, our proofs rely on the analysis of the asymptotic behavior of the cumulants of the random variables concerned  (recall that the distribution of the elements of the second Wiener is completely determined by their cumulants, see \cite{FT} or \cite{NPbook}).

We organized our paper as follows. In Section 2 we recall some basic definitions for the Rosenblatt process and the Wiener-Rosenblatt integral and we also state a general result for the limit behavior in law of the Wiener-Rosenblatt integral  as the Hurst index converges to its extreme values. In Section 3 we treat the particular case of the Rosenblatt Ornstein-Uhlenbeck process.

\section{The Rosenblatt process and the Wiener-Rosenblatt integral}

Below, in the first part,  we present the definition and the basic properties of the Rosenblatt process and of the Wiener integral with respect to the Rosenblatt process. In the second part, we give a general result concerning the convergence in distribution with respect to the Hurst parameter of the Wiener-Rosenblatt integral. This general result will be applied in the next section in order to obtain the limit behavior of the Ornstein-Uhlenbeck process with Rosenblatt noise.

\subsection{Definition and basic properties}
Let $(\Omega, \mathcal{F}, P)$ be a probability space and let $(B(y))_{y \in \mathbb{R}} $ be a Wiener process on $\Omega$. We will denote by $(Z^{H}(t)) _{t\geq 0} $ the Rosenblatt process with self-similarity index $H\in \left( \frac{1}{2}, 1\right)$. It is defined on $\Omega$ as a multiple stochastic integral of order 2 with respect to the Wiener process $B$ via
\begin{equation}
\label{rose1}
Z ^{H}(t)= c(H, 2) \int_{\mathbb{R}} \int_{\mathbb{R}} \left(    \int_{0} ^{t} (u-y_{1})_{+} ^{\frac{H}{2}-1}(u-y_{2})_{+} ^{\frac{H}{2}-1}du\right)dB(y_{1}) dB(y_{2})= I_{2} (L ^{H} _{t}), \hskip0.3cm t\geq 0
\end{equation}
where $I_{2}$ denotes the multiple stochastic integral of order two with respect to $B$ (see the Appendix) and  we denoted by $ L ^{H}_{t}$ the kernel of the Rosenblatt process given by, for every $y_{1}, y_{2}\in \mathbb{R}$
\begin{equation}
\label{lh}
L ^{H}_{t}(y_{1}, y_{2})= c(H,2)\int_{0} ^{t} (u-y_{1})_{+} ^{\frac{H}{2}-1}(u-y_{2})_{+} ^{\frac{H}{2}-1}du.
\end{equation}
We denoted $x_{+}=\max (x, 0)$. It is well-known (see e.g. \cite{T}) that the kernel $ L ^{H} _{t} $ belongs to $ L ^{2} (\mathbb{R} ^{2})$ for every $t\geq 0$ when $H>\frac{1}{2}$,  which implies that the multiple integral of order two in (\ref{rose1}) is well-defined. The strictly positive constant $c(H,2)$ is chosen such that $ \mathbf{E} (Z ^{H}(1)) ^{2}=1$. Actually  (see e.g. \cite{T}, Proposition 3.1)

\begin{equation}
\label{d2}
c( { H}, 2) ^{2}=  \frac{ {H} (2{
 H}-1) } {2 \beta \left( \frac{ {  H}}{2}, 1-{H}\right)^{2}}, 
\end{equation}
where $\beta $ is Beta function $\beta (p,q)= \int_{0}^{1} z^{p-1}(1-z)^{q-1} dz, p,q>0$. 

The Rosenblatt process is a $H$-self-similar process with stationary increments. It exhibits long-range dependence and its sample paths are H\"older continuous of order $\delta$ for any $\delta \in (0, H)$.   This process has been intensively studied in the last decades, see e.g. the monographs \cite{PiTa-book} and \cite{T} and the references therein.

The Wiener integral with respect to the Rosenblatt process (or the {\em Wiener-Rosenblatt integral}) has been constructed in \cite{MaTu}.  Let $f \in \mathcal{H}_{H}$ where $\mathcal{H}_{H}=\{ f: \mathbb{R}\to \mathbb{R}: \Vert f\Vert _{\mathcal{H}_{H}} <\infty \}  $  with
\begin{equation*}
\Vert f\Vert _{\mathcal{H}_{H}} ^{2}: =H(2H-1) \int_{\mathbb{R}}\int_{\mathbb{R}} f(u)f(v) \vert u-v\vert ^{2H-2} dudv.
\end{equation*}
Then the Wiener integral of $f$ with respect to $Z ^{H}$ is given by 
\begin{equation}
\label{wr}
\int_{\mathbb{R}} f(u) d Z ^{H}(u)= I_{2} (J_{H}f)
\end{equation}
where the kernel $J_{H}f $ has the expression, for every $y_{1}, y_{2} \in \mathbb{R}$, 
\begin{equation}
\label{jf}
(J_{H}f) (y_{1}, y_{2}) =c(H, 2) \int_{\mathbb{R}} f(u) (u-y_{1})_{+} ^{\frac{H}{2}-1}(u-y_{2})_{+} ^{\frac{H}{2}-1}du.
\end{equation}
The Wiener-Rosenblatt integral constitutes an isometry between  $\mathcal{H}_{H}$ and $L^{2} (\Omega)$ since, for every $f,g\in \mathcal{H}_{H}$
\begin{equation}
\label{isow}
\mathbf{E}\left(  \int_{\mathbb{R}} f(u) d Z ^{H}(u)\int_{\mathbb{R}} g(u) d Z ^{H}(u)\right)= \langle f,g\rangle _{\mathcal{H}_{H}} : = H(2H-1)  \int_{\mathbb{R}}\int_{\mathbb{R}} f(u)g(v) \vert u-v\vert ^{2H-2} dudv.
\end{equation}
The space $\mathcal{H} _{H}$ is not complete and may contains distributions. A subspace of functions included in $\mathcal{H}_{H}$ is the space $\left| \mathcal{H}_{ H}\right| $ of measurable functions $f: \mathbb{R} \to \mathbb{R}$ such that 
\begin{equation}
\label{hh2}
\Vert f\Vert ^{2}_{\left|  \mathcal{H}_{H} \right|  } := \int_{\mathbb{R}}\int_{\mathbb{R}} \vert  f(u)f(v)\vert  \vert u-v\vert ^{2H-2} dudv<\infty.
\end{equation}
\subsection{Asymptotic behavior of the Wiener-Rosenblatt  integral}

Our purpose is to study the asymptotic behavior, as $H\to \frac{1}{2}$ and $H\to 1$, of the Wiener integral with respect ot the Rosenblatt process$$\int_{\mathbb{R}}f(u) dZ ^ {H}(u) $$
where
$(Z ^ {H}(t))_{t\geq 0}$ is a Rosenblatt process with self-similarity order $H\in (\frac{1}{2}, 1)$ and $f\in \mathcal{H}_{H}$.

The proof of the asymptotic behavior of the Wiener-Rosenblatt integral is based on the analysis of its cumulants.  Since the random variable $$\int_{\mathbb{R}}f(u) dZ ^ {H}(u) $$ belongs to the second Wiener chaos, see (\ref{wr}), its law is completely determined by its cumulants (or equivalently, by its moments). That is, if $F,G$ are elements of the second Wiener chaos then $F$ and $G$ have the same law if and only if they have the same cumulants. Moreover, the convergence of the cumulants implies the convergence in distribution when we deal with sequences  in the  Wiener chaos of order two. Let us denote by $k_{m} (F)$, $m\geq 1$ the  $m$th cumulant of a random variable $F$. It is defined as
\begin{equation*}
k_{m}(F)= (-i) ^{n} \frac{\partial ^{n}}{\partial t^{n} }\ln \mathbf{E} (e ^{itF}) | _{t=0}.
\end{equation*}
We have the following link between the moments and the cumulants of $F$: for every $m\geq 1$, 
\begin{equation}\label{mom-cum}
k_{m}(F)= \sum_{\sigma =(a_{1},.., a_{r})\in \mathcal{P}( \{1,..,n\} )}(-1) ^ {r-1} (r-1)! \mathbf{E} X ^ {\vert a_{1}\vert}\ldots \mathbf{E} X ^ {\vert a_{r}\vert}
\end{equation}
if $F\in L^{m}(\Omega)$, where $\mathcal{P} (b)$ is the set of all partitions of $b$. In particular, for centered random variables $F$, we have $k_{1}(F)=\mathbf{E} F, k_{2}(F)= \mathbf{E} F^{2},  k_{3}(F)= \mathbf{E}F^{3},  k_{4}= \mathbf{E} F^{4}- (\mathbf{E} F^{2}) ^{2}.$

In the particular situation when $G=I_{2}(f) $ is a multiple integral of order 2 with respect to a Wiener process  $(B(y))_{y\in \mathbb{R}  }$, then its cumulants can be computed as (see e.g. \cite{Nourdin}, Proposition 7.2 or \cite{T})
\begin{equation}
k_{m}(G)=  2^{m-1}(m-1)! \int_{ \mathbb{R}  ^{m} }d { u}_{1}\ldots d{ u}_{m} f({ u}_{1},
{ u}_{2}) f({ u}_{2}, { u}_{3})\ldots f({ u}_{m-1}, { u}_{m}) f({ u}_{m},{  u}_{1}).
\label{cum2}
\end{equation}

From the formula (\ref{cum2}), we can obtain the following expression of the cumulants of the Wiener-Rosenblatt integral (see e.g. \cite{ST}).

\begin{prop}We have, for $m\geq 2$
\begin{eqnarray}
&&k_{m} \left( \int_{\mathbb{R}} f(u) d Z^ {H}(u) \right) \label{cum1}\\
&&= c_{1,m} \int_{\mathbb{R}}\ldots \int_{\mathbb{R}} du_{1}...du_{m} f(u_{1})...f(u_{m}) \vert u_{1}- u_{2} \vert ^ {H-1} \vert u_{2}-u_{3} \vert ^ {H-1}...\vert u_{m-1}-u_{m}\vert ^{H-1} \vert u_{m}-u_{1}\vert ^ {H-1}\nonumber
\end{eqnarray}
with 
\begin{equation}\label{c1m}
c_{1,m}= 2 ^ {\frac{m}{2}-1}(m-1) ! (H(2H-1))^ {\frac{m}{2}}
\end{equation}
and $k_{1} \left( \int_{\mathbb{R}} f(u) d Z^ {H}(u) \right) =0.$

\end{prop}

We will treat separately the behavior of the Wiener-Rosenblatt integral as $H$ is  near $\frac{1}{2}$ and near $1$. The limiting process will be different in these two cases.

\subsubsection{Convergence when $H\to 1$}

We have the following result.

\begin{prop}\label{p2}
Consider $f:\mathbb{R} \to \mathbb{R} $  such that for some $\eps \in (0, \frac{1}{2})$

\begin{equation}\label{cH1}
\Vert f\Vert ^{2}  _{\left|   \mathcal{H} _{\frac{1}{2}+\eps}\right| }= \int_{\mathbb{R}}\int_{\mathbb{R}} dudv\vert f(u) \vert \vert f(v) \vert \vert u-v \vert ^{2\eps -1} <\infty \mbox{ and } f\in L^ {1}(\mathbb{R}).
\end{equation}
Then

\begin{equation}\label{8m-1}
 \int_{\mathbb{R}} f(u) d Z^ {H}(u)\xrightarrow[H\to 1]{(d)} \frac{1}{\sqrt{2}}\left(  \int_{\mathbb{R}} f(u)du\right) (Z ^ {2}-1)
\end{equation}
with $Z\sim N(0,1)$. 
\end{prop}
{\bf Proof: }  First notice that condition (\ref{cH1}) implies that $f\in \left| \mathcal{H} _{ H }\right|$ for $H\geq \frac{1}{2}+\eps$. Indeed, by using the bound $\sup_{H\in [\frac{1}{2}+\eps, 1]} \vert x\vert ^{2H-2} \leq 1\vee \vert x\vert ^{2\eps -1} $, we get 
\begin{eqnarray*}
\Vert f\Vert ^{2} _{\vert \mathcal{H} _{H}\vert }&\leq & \int_{\mathbb{R}} \int_{\mathbb{R}} dudv \vert f(u)f(v)\vert (1\vee \vert u-v\vert ^{2\eps -1}) \\
&\leq & \left( \int_{\mathbb{R} }\vert f(u)\vert du\right) ^{2} + \int_{\mathbb{R}}\int_{\mathbb{R}} dudv\vert f(u) \vert \vert f(v) \vert \vert u-v \vert ^{2\eps -1} <\infty.
\end{eqnarray*}

Consider the random variable 
$$G:= I_{2} \left(  \frac{1}{\sqrt{2}}\left(  \int_{\mathbb{R}} f(u)du\right) 1_{[0,1]} ^ {\otimes 2}\right)$$
which has the same law as the right-hand side of (\ref{8m-1}). We have
$$k_{m}(G)= 2 ^ {m-1}(m-1)!\left( \frac{1}{\sqrt{2}}\left(  \int_{\mathbb{R}} f(u)du\right) \right) ^ {m}= 2 ^ {\frac{m}{2}-1}(m-1) ! \left(  \int_{\mathbb{R}} f(u)du \right) ^ {m} .$$
On the other hand, since by (\ref{c1m})
$$c_{1,m }\xrightarrow[H\to 1]{}  2 ^ {\frac{m}{2}-1}(m-1) ! $$
the conclusion will follow if we show that 

$$ \int_{\mathbb{R}}...\int_{\mathbb{R}} du_{1}...du_{m} f(u_{1})...f(u_{m}) \vert u_{1}- u_{2} \vert ^ {H-1} \vert u_{2}-u_{3} \vert ^ {H-1}...\vert u_{m-1}-u_{m}\vert ^{H-1} \vert u_{m}-u_{1}\vert ^ {H-1}$$
converges, as $H\to 1$, to 
$$\left( \int_{\mathbb{R} }f(u) du \right) ^ {m}. $$

Assume that the integrability condition (\ref{cH1})  is satisfied.  Consider  the function $g_{H}$ defined on   $\mathbb{R} ^{m}$ except the diagonals with values in $ \mathbb{R}$, given by 
$$g_{H}(u_{1},.., u_{m}) = f(u_{1})...f(u_{m}) \vert u_{1}- u_{2} \vert ^ {H-1} \vert u_{2}-u_{3} \vert ^ {H-1}...\vert u_{m-1}-u_{m}\vert ^{H-1} \vert u_{m}-u_{1}\vert ^ {H-1}.$$
Clearly $g_{H}(u_{1},.., u_{m}) $ converges to $ f(u_{1})...f(u_{m})$ for almost  every $u_{1},.., u_{m}\in \mathbb{R}$.
Also, using again the bound $\sup_{H\in [\frac{1}{2}+\eps, 1]} \vert x \vert ^{H-1} \leq 1\vee \vert x\vert ^{-\frac{1}{2}+\eps} $ for every $x\in \mathbb{R}$, we have
 
\begin{eqnarray*}
\sup_{H\in [\frac{1}{2}+\eps, 1]}\vert g_{H} (u_{1},.., u_{m}) \vert 
&&\leq  \vert f(u_{1})...f(u_{m}) \vert \left( 1\vee (\vert u_{1}-u_{2}\vert ^ {-\frac{1}{2}+\eps}.... \vert u_{m}-u_{1}\vert ^ {-\frac{1}{2}+\eps} )\right) \\
&&\leq  \vert f(u_{1})...f(u_{m}) \vert +  \vert f(u_{1})...f(u_{m}) \vert.....\vert u_{1}-u_{2}\vert ^ {-\frac{1}{2}+\eps}.... \vert u_{m}-u_{1}\vert ^ {-\frac{1}{2}+\eps}.
\end{eqnarray*}
In order to apply the dominated convergence theorem, we need to show that the function
$$\vert f(u_{1})...f(u_{m}) \vert +  \vert f(u_{1})...f(u_{m}) \vert.....\vert u_{1}-u_{2}\vert ^ {-\frac{1}{2}+\eps}.... \vert u_{m}-u_{1}\vert ^ {-\frac{1}{2}+\eps}$$
is integrable over $\mathbb{R} ^{m}$. Since $f\in L ^{1} (\mathbb{R}) $, the first summand  above is integrable over $\mathbb{R} ^{m}$. On the other hand, $\int_{\mathbb{R}}...\int_{\mathbb{R} } du_{1}..du_{m}  \vert f(u_{1})...f(u_{m}) \vert.....\vert u_{1}-u_{2}\vert ^ {-\frac{1}{2}+\eps}.... \vert u_{m}-u_{1}\vert ^ {-\frac{1}{2}+\eps}$ represents, modulo a constant, the $m$th cumulant of the random variable $I_{2} (J_{\frac{1}{2}+\eps} \vert f\vert)$, where $ J_{\frac{1}{2}+\eps} $ is given by (\ref{jf}). The fact that $\Vert f\Vert ^{2}  _{\left|   \mathcal{H} _{\frac{1}{2}+\eps}\right| }<\infty $ (see (\ref{cum1})) by (\ref{cH1}) and the hypercontractivity property of multiple integrals (\ref{hyper})  imply that the all the moments of $I_{2} (J_{\frac{1}{2}+\eps} \vert f\vert)$ are finite and consequently, via (\ref{mom-cum}), all the cumulants of $I_{2} (J_{\frac{1}{2}+\eps} \vert f\vert)$ are finite. This implies the conclusion by using the dominated convergence theorem. \qed

\begin{remark}
For example, when $f$ is bounded with compact support, then condition (\ref{cH1}) is satisfied (in particular if $f=1_{[0, t]}$ with $t>0$ fixed).  A detailed explanation  can be found in the proof of Proposition \ref{p4} below.
\end{remark}

\subsubsection{Convergence when $H\to \frac{1}{2}$ }

Concerning the behavior of the Rosenblatt-Wiener integral when $H$ approaches one half, we have the following result.

\begin{prop}\label{p3}
Assume $f\in\mathcal{H}_{H}$. Also assume that
\begin{equation}
\label{sf}
\sigma _{f} ^{2}= \lim _{H\to \frac{1}{2} } \Vert f\Vert ^{2} _{\mathcal{H} _{H}}=\lim _{H\to \frac{1}{2} }  H(2H-1) \int_{\mathbb{R}}\int_{\mathbb{R}} f(u)f(v) \vert u -v\vert ^{2H-2}dudv
\end{equation}
exists and it is finite and
\begin{equation}\label{cH12}
(2H-1) ^{2}   \int_{\mathbb{R}^{4}} du_{1}...du_{4} f(u_{1})...f(u_{4}) \vert u_{1}- u_{2} \vert ^ {H-1} \vert u_{2}-u_{3} \vert ^ {H-1}\vert u_{3}- u_{4} \vert ^ {H-1}  \vert u_{4}-u_{1}\vert ^ {H-1} \xrightarrow[H\to  \frac{1}{2}]{}0.
\end{equation}
Then
\begin{equation*}
\int_{\mathbb{R}}f(u) d Z^ {H}(u) \xrightarrow[H\to  \frac{1}{2}]{(d)}  N (0, \sigma_{f} ^{2}).
\end{equation*}

\end{prop}
{\bf Proof: } We need to analyze the behavior of the cumulants of  $\int_{\mathbb{R}}f(u) d Z^ {H}(u)$ as $H$ converges to $\frac{1}{2}$.  Recall that $k_{1} \left( \int_{\mathbb{R}}f(u) d Z^ {H}(u)\right)= \mathbf{E}\left( \int_{\mathbb{R}}f(u) d Z^ {H}(u)\right)=0$ and
$$k_{2} \left( \int_{\mathbb{R}}f(u) d Z^ {H}(u)\right) = H(2H-1) \int_{\mathbb{R}} \int_{\mathbb{R}} du_{1}du_{2} f(u_{1}) f(u_{2} )\vert u_{1}-u_{2}\vert ^ {2H-2}=\Vert f\Vert _{\mathcal{H}_{H}} ^ {2}. $$

By (\ref{sf}), we have $k_{2} \left( \int_{\mathbb{R}}f(u) d Z^ {H}(u)\right)\xrightarrow[H\to  \frac{1}{2}]{}\sigma _{f} ^{2}$.  On the other hand, due to condition (\ref{cH12}), 
$$k_{4}  \left( \int_{\mathbb{R}}f(u) d Z^ {H}(u)\right) \xrightarrow[H\to  \frac{1}{2}]{} 0.$$
Since $\int_{\mathbb{R}}f(u) d Z^ {H}(u)$ belongs to the second Wiener chaos,  by the Fourth Moment Theorem (see \cite{NuPe}, see also Theorem \ref{fmt} in the Appendix) we conclude that $\int_{\mathbb{R}}f(u) d Z^ {H}(u)$ converges in distribution as $H\to \frac{1}{2}$ to a Gaussian random variable with mean zero (the limit of its first cumulant) and variance $\sigma_{f} ^{2}$ (the limit of its second cumulant).  \qed

\vskip0.3cm

We will see below that for certain kernels $f$ the condition (\ref{sf}) is automatically satisfied.  Recall that a sequence of functions $(f_{n})_{n\geq 1}$ is an approximation of the identity as $n\to \infty$ if \begin{itemize}
\item $f_{n}(t)\geq 0$ for every $t\in \mathbb{R}$.
\item For every $\delta >0$, $\int_{\vert t\vert \leq \delta } f_{n}(t)dt \xrightarrow[n\to \infty]{}1. $
\item For every $\delta >0$, $\int_{\vert t\vert > \delta } f_{n}(t)dt \xrightarrow[n\to \infty]{}0. $
\end{itemize}
Moreover, if $(f_{n})_{n\geq 1}$ is an approximation of the identity, then  for every $f\in L ^{p} (\mathbb{R})$  with $p\in [1, \infty)$, the convolution $f\ast f_{n}$ converges in $ L^{p} (\mathbb{R} ) $ to $f$.

\begin{corollary}
\label{cor1}
Assume $f\in\mathcal{H}_{H}\cap L ^ {2} (\mathbb{R})$ with $\supp (f) \subset [0, \infty)$ and (\ref{cH12}) holds. Then 
\begin{equation*}
\int_{\mathbb{R}}f(u) d Z^ {H}(u) \xrightarrow[H\to  \frac{1}{2}]{(d)} N (0, \int_{\mathbb{R}} f^{2}(u) du).
\end{equation*}
\end{corollary}
{\bf Proof: } By Proposition \ref{p2}, it suffices to check that the limit (\ref{sf}) exists and it is equal to $\int_{\mathbb{R}} f^{2}(u) du$. We have 
\begin{eqnarray*}
\Vert f\Vert _{\mathcal{H}_{H}}^{2}&=& H(2H-1) \int_{0} ^{\infty} \int_{0} ^{\infty} dudv f(u)f(v)\vert u-v\vert ^{2H-2} \\
&=& 2H (2H-1) \int_{0}^{\infty} du f(u) \int_{0} ^{u}  dv f(u-v) v ^{2H-2}. 
\end{eqnarray*}
Notice that the function $2H(2H-1)1_{[0,u]}(v) v ^{2H-2}$ constitutes an approximation of the identity as $H\to \frac{1}{2}$. Therefore,  
$$\Vert f\Vert _{\mathcal{H}_{H}}^{2} \xrightarrow[H\to  \frac{1}{2}]{} \int_{\mathbb{R}} f^{2}(u)du.$$
By condition (\ref{cH12}) and Theorem \ref{fmt},  the conclusion follows.  \qed

\begin{remark}
\begin{itemize}
\item
The above result shows that, when $H\to \frac{1}{2}$, the Wiener-Rosenblatt integral $\int_{\mathbb{R}}f(u) d Z^ {H}(u)$ converges in distribution to $\int_{\mathbb{R}} f(u) dW(u)$, where $W$ is a Wiener process. This is a natural extension of the results in \cite{ArTu} or \cite{VeTa}.

\item The condition $\supp(f)\subset [0, \infty)$  in Corollary \ref{cor1} cannot be omitted. If for example $\supp (f)$ is $\mathbb{R}$, the above proof does not work, since $2H(2H-1)1_{[0,\infty)}(v) v ^{2H-2}$ is not an approximation of the unity.
\end{itemize} 
\end{remark}

\section{Asymptotic behavior of the Rosenblatt Ornstein-Uhlenbeck process}
The Rosenblatt Ornstein-Uhlenbeck (ROU)  process is defined as the unique solution of the Langevin equation 

\begin{equation}
\label{ou1}
X _{t}= \xi -\lambda \int_{0} ^ {t} X_{s}ds + \sigma Z ^ {H}(t), \hskip0.4cm t\geq 0
\end{equation}
where $\lambda, \sigma >0$ and the initial condition $\xi$ is a random variable in $ L ^ {2} (\Omega)$. The case when the noise in (\ref{ou1}) is the fractional Brownian motion has been considered in \cite{CKM}. 
 
The unique solution to (\ref{ou1}) can be expressed as 
\begin{equation}
\label{ou2}
Y ^ {H}(t) =e ^ {-\lambda t}\left( \xi + \sigma \int_{0} ^ {t} e ^ {\lambda u} d Z^ {H}(u)\right)
\end{equation}
where the stochastic integral with respect to $ Z ^ {H}$ can be understood both in the Wiener or Riemann-Stieltjes sense.

The stationary Rosenblatt Ornstein-Uhlenbeck process is obtained by taking the initial condition $\xi= \sigma \int_{-\infty} ^{0} e ^{-\lambda u} d Z ^{H}(u) $ in (\ref{ou1}). Then, the stationary ROU, which will be denoted in the sequel by $ (X ^{H}(t))_{t\geq 0}$, can be expressed as, for every $t\geq 0$,
\begin{equation}
\label{ou3}
X ^{H}(t)= \sigma\int_{-\infty} ^{t}  e ^{-\lambda (t-u) }d Z ^{H}(u).
\end{equation}
The process $(X ^{H}(t))_{t\geq 0}$ is a stationary Gaussian process, $H$-self-similar with stationary increments. Moreover, it exhibits long-range dependence since $H>\frac{1}{2}$, see \cite{CKM} or \cite{MaTu}.

In this paragraph, our purpose is to analyze the asymptotic behavior, as $ H\to 1$ and as $H\to \frac{1}{2}$, in the sense of the weak convergence,  of the processes (\ref{ou2}) and (\ref{ou3}). The analysis of the limit behavior at the extreme critical values of the Hurst exponent is different for $X ^{H}$ and $Y^{H}$ due to the fact that the the initial values depends on $H$ in the case of $ X ^{H}$.

\subsection{Padded sets and the power counting theorem}\label{power}

We need to recall some notation and results from \cite{TeTa} which are needed in order to check the integrability assumption from Proposition  \ref{p3}.

Consider a set $T=\{M_{1},.., M_{m}\}$ of linear functions on $\mathbb{R } ^ {m}$.  The power counting theorem (see Theorem 1.1 and Corollary 1.1 in \cite{TeTa}) gives sufficient conditions for the integral 
\begin{equation}
\label{i}
I= \int_{\mathbb{R}}...\int_{\mathbb{R}} du_{1}...du_{m} f_{1} (M_{1} (u_{1},.., u_{m}))....f_{m} (M_{m} (u_{1},.., u_{m}))
\end{equation}
to be finite, where $f_{i}:\mathbb{R} \to \mathbb{R}$, $i=1,.., m$ are such that $\vert f_{i} \vert $ is bounded above on $(a_{i}, b_{i})$ ($0<a_{i}<b_{i} <\infty$) and
$$ \vert f_{i}(y) \vert \leq c_{i} \vert y\vert ^ {\alpha _{i} } \mbox{ if } \vert y_{i}\vert <a_{i} \mbox{ and } \vert f_{i}(y)\vert \leq c_{i}\vert y\vert ^ {\beta _{i} } \mbox{ if } \vert y\vert >b_{i}.$$

For a subset $W\subset T$ we denote  by $s_{T}(W) =span (W)\cap T$. A subset $W$ of $T$ is said to be {\em padded } if $s_{T}(W)=W$ and any functional $M\in W$ also belongs to $s_{T}(W\setminus \{M\}).$ Denote by $\span (W)$ the linear span generated by $W$ and by $r(W)$ the number of linearly independent elements of $W$. 

Then Theorem 1.1  in \cite{TeTa} says that the integral $I$ (\ref{i}) is finite if
\begin{equation}
\label{d0}
 d_{0} (W)= r(W) + \sum _{s_{T}(W)} \alpha _{i} >0 
\end{equation}
for any  subset $W$ of $T$ with   $s_{T}(W)=W$ and 
\begin{equation}
\label{dinf}
 d_{\infty } (W)= r(T) -r(W) + \sum _{T\setminus s_{T}(W)} \beta _{i} <0
\end{equation}
for any proper subset $W$ of $T$ with $s_{T}(W)= W$, including the empty set. If $\alpha _{i}>-1$ then  it suffices to check (\ref{d0}) for any padded subset $W\subset T$. Also, it suffices to verify (\ref{dinf}) only for padded subsets of $T$ if $ \beta _{i} \geq -1.$

The condition (\ref{d0}) implies the integrability at the origin while (\ref{dinf}) gives the integrability of $I$ at infinity.

There is a similar result if one starts with a set $T$ of affine functionals instead of linear functionals. 

\subsection{ The (non-stationary) ROU process}
We first treat the case of the process (\ref{ou2}) with initial condition not depending on the Hurst index. In the sequel, we fix $T>0$ arbitrary chosen.

\subsubsection{Convergence when $H\to 1$}

\begin{prop}\label{p4}
Assume that the initial condition $\xi$  does not depend on $H$. Then the process $(Y ^ {H}(t))_{t\in [0,T]} $ converges weakly, in the space of continuous functions $C[0,T]$, to the stochastic process $(Y(t))_{t\in [0,T]}$ given by 
\begin{equation}
\label{y}
 Y(t)= e ^ {-\lambda t}\xi + \sigma \frac{1}{\sqrt{2}} \left( \int_{0}^ {t} e ^ {-\lambda (t-u) } du \right)(Z ^ {2}-1)=  e ^ {-\lambda t}\xi +\frac{\sigma }{\sqrt{2}\lambda} (1-e ^ {-\lambda t} )(Z ^ {2}-1) 
\end{equation}
with $Z \sim N(0,1).$
\end{prop}
{\bf Proof: } We start by checking the convergence of the finite dimensional distribution of $Y ^{H}$ of those of $Y$. Take  $\alpha_{1}, ..., \alpha _{d}\in \mathbb{R}$ and $t_{1},..., t_{d} \in [0, T]$. We will prove that 
$\sum_{j=1}^ {d} \alpha _{j} Y ^ {H}(t_{j}) $ converges in distribution, as $H\to 1$ to the random variable $\sum_{j=1}^ {d} \alpha _{j} Y (t_{j})$. 

We have, by the linearity of the Wiener-Rosenblatt integral,
$$\sum_{j=1}^ {d} \alpha _{j} Y ^ {H}(t_{j}) = \int_{\mathbb{R}} \left( \sum_{j=1}^ {d} \alpha _{j} 1_{[0, t_{j}]} (u) e ^{ -\lambda (t_{j}- u)}\right) dZ ^ {H}(u)=\int_{\mathbb{R}}f(u)d Z^{H}(u)$$
with
\begin{equation}
\label{f1}
f(u)= \sum_{j=1}^ {d} \alpha _{j} 1_{[0, t_{j}]} (u) e ^{ -\lambda (t_{j}- u)}.
\end{equation}
In order to apply Proposition \ref{p2}, we need to show that condition (\ref{cH1})  holds true.  Clearly $f$ belongs to $L ^ {1} (\mathbb{R})$.  Concerning the first part of (\ref{cH1}), we have
\begin{eqnarray*}
&& \int_{\mathbb{R}}\int_{\mathbb{R}}dudv \vert f(u)f(v)\vert \vert u-v\vert ^{2\eps-1} \leq \sum_{j,k=1} ^{d} \vert \alpha _{j} \alpha _{k}\vert \int_{0} ^{t_{j}} du \int_{0} ^{t_{k}} dv e ^{-\lambda (t_{j}-u)} e ^{-\lambda (t_{k}-v)} \vert u-v\vert ^{2\eps-1}\\
&\leq &\sum_{j,k=1} ^{d} \vert \alpha _{j} \alpha _{k}\vert \int_{0} ^{T} du \int_{0} ^{T} dv e ^{-\lambda (t_{j}-u)} e ^{-\lambda (t_{k}-v)} \vert u-v\vert ^{2\eps-1}\leq \sum_{j,k=1} ^{d} \vert \alpha _{j} \alpha _{k}\vert \int_{0} ^{T} du \int_{0} ^{T} dv \vert u-v\vert ^{2\eps-1}\\
&=&\frac{1}{\eps(\eps+1)}T ^{2\eps +1}  \sum_{j,k=1} ^{d} \vert \alpha _{j} \alpha _{k}\vert<\infty.
\end{eqnarray*}

Concerning the tightness, notice that for every $0\leq s<t\leq T$ we have, since $Y^{H}$ is a solution to (\ref{ou1})
$$\mathbf{E} \vert Y^ {H}(t)- Y ^ {H}(s)\vert ^ {2} \leq C\vert t-s\vert $$
and therefore, for every $p\geq 1$,
\begin{equation}\label{7d-1}
\mathbf{E}\vert Y^ {H}(t)- Y^ {H}(s) \vert ^ {2p}\leq C_{p} (\mathbf{E} \vert Y^ {H}(t)- Y ^ {H}(s)\vert ^ {2}) ^ {p}\leq c\vert t-s\vert ^ {p}.
\end{equation}
The tightness is obtained from (\ref{7d-1}) and e.g. Lemma 2.2 in \cite{Pro}. \qed

 \vskip0.2cm

Note that the limit process $(Y(t))_{t\in [0, T]} $ given by (\ref{y}) is a second chaos stochastic process. Therefore, its finite dimensional distributions are characterized by the cumulants, i.e. for every $\alpha_{1},.., \alpha _{d}\in \mathbb{R}$ and $t_{1},.., t_{d}\in [0, T]$, 
\begin{equation*}
k_{m}\left(   \sum_{j=1} ^{d} \alpha _{j} Y (t_{j}) \right)=k_{m} \left( \frac{\sigma }{\sqrt{2} \lambda }1_{[0,1] }^{\otimes 2}  \sum_{j=1} ^{d} \alpha _{j} (1- e ^{-\lambda t_{j}}) \right)=2 ^{\frac{m}{2}-1} (m-1) ! \frac{ \sigma ^{m} }{ \lambda ^{m}} \left( \sum_{j=1} ^{d} \alpha _{j} (1- e ^{-\lambda t_{j}}) \right) ^{m},
\end{equation*}
for $m\geq 2$ and 
$$k_{1}\left(   \sum_{j=1} ^{d} \alpha _{j} Y (t_{j}) \right)= \mathbf{E}\left(   \sum_{j=1} ^{d} \alpha _{j} Y (t_{j}) \right)=\mathbf{E}\left( \xi \right). \sum_{j=1} ^{d} \alpha _{j} e ^{-\lambda t_{j}} .$$

\subsubsection{Convergence when $H\to \frac{1}{2}$}

The (standard) Ornstein-Uhlenbeck process (denoted $Y_{0}$ in the sequel) is given by (\ref{ou2}) with $Z ^{H}$ replaced by a Wiener process $W$. Thus it can be written as 
\begin{equation}\label{y0}
Y_{0}(t)=e ^{-\lambda t} \left( \xi +\sigma \int_{0} ^{t} e ^{\lambda u} dW(u) \right), \mbox{ for every }t\geq 0.
\end{equation}
Consequently, $(Y_{0}(t))_{t\geq 0} $ is a Gaussian process with mean $\mathbf{E} Y_{0}(t)= e ^{-\lambda t} \mathbf{E}\xi$ for any $t\geq 0$ and covariance function
\begin{equation*}
Cov (Y_{0}(t), Y_{0} (s))= \frac{\sigma ^{2}} {2\lambda } \left( e ^{-\lambda \vert t-s\vert}- e ^{-\lambda (t+s)}\right)
\end{equation*}
for every $s,t\geq 0$.

The  Ornstein-Uhlenbeck process will appear as limit of the  ROU  process as $H\to \frac{1}{2}$.

\begin{prop}
As $H\to \frac{1}{2}$, the process $Y^ {H}(t) $ converges weakly to the Ornstein-Uhlenbeck process $(Y_{0} (t))_{t\in [0, T]}$. 
\end{prop}
{\bf Proof: } Consider $\lambda _{1},.., \lambda _{d} \in \mathbb{R}$ and $t_{1},.., t_{d}\in [0, T] $.  We will apply  Corollary \ref{cor1}. Clearly, as in the proof of Proposition \ref{p4},  the function $f(u)= \sum_{j=1}^ {d} \alpha _{j} 1_{[0, t_{j}]} (u) e ^{ -\lambda (t_{j}- u)}$ from (\ref{f1}) belongs to $\mathcal{H}_{H}\cap L ^ {2} (\mathbb{R})$.

We need to show that condition (\ref{cH12}) is satisfied. We have

\begin{eqnarray*}
&&  \int_{\mathbb{R}^{4}} du_{1}...du_{4} f(u_{1})...f(u_{4}) \vert u_{1}- u_{2} \vert ^ {H-1} \vert u_{2}-u_{3} \vert ^ {H-1}\vert u_{3}-u_{4}\vert ^{H-1} \vert u_{4}-u_{1}\vert ^ {H-1}\\
&\leq &\sum_{j_{1},.., j_{4}=1}^ {d} \vert \alpha _{j_{1}}....\alpha_{j_{4}}\vert  \int_{0} ^ {T}...\int_{0} ^ {T} du_{1}..du_{4}  \vert u_{1}- u_{2} \vert ^ {H-1} \vert u_{2}-u_{3} \vert ^ {H-1}\vert u_{3}-u_{4}\vert ^{H-1} \vert u_{4}-u_{1}\vert ^ {H-1}.
\end{eqnarray*}
Actually, the function 
$$H\to  \int_{0} ^ {T}...\int_{0} ^ {T} du_{1}...du_{4}  \vert u_{1}- u_{2} \vert ^ {H-1} \vert u_{2}-u_{3} \vert ^ {H-1}\vert u_{3}-u_{4}\vert ^{H-1} \vert u_{4}-u_{1}\vert ^ {H-1}$$ is finite and continuous on the set $(\frac{1}{4}, 1]$. This follows from Lemma 3.3 in \cite{BaiTa} but also by applying the power counting theorem with $(\alpha_{1},.., \alpha_{4})=(H-1,.., H-1)$. Therefore
$$\sup_{H\in [\frac{1}{2}, 1]}   \int_{\mathbb{R}}...\int_{\mathbb{R}} du_{1}...du_{4} f(u_{1})...f(u_{4}) \vert u_{1}- u_{2} \vert ^ {H-1} \vert u_{2}-u_{3} \vert ^ {H-1}\vert u_{3}-u_{4}\vert ^{H-1} \vert u_{4}-u_{1}\vert ^ {H-1}<\infty$$
and  implies (\ref{cH12}) and consequently, it gives the convergence of finite dimensional distributions of $Y^{H}$ as $ \to \frac{1}{2}$. The tighness follows  (\ref{7d-1}). \qed

\subsection{Asymptotic behavior of the stationary Rosenblatt Ornstein-Uhlenbeck process}

Now, we analyze the asymptotic behavior of the process (\ref{ou3}). The main idea is the same as in the previous section but we need to pay attention to the fact that the kernel of the stationary ROU process has the whole real line as support .

\subsubsection{Convergence when $H\to 1$}

\begin{prop}
The stationary  Rosenblatt Ornstein-Uhlenbeck process converges weakly, in the space of continuous functions $C[0, T]$, to the stochastic process $ (X(t))_{t\in [0, T]}$ defined by, for every $t\in [0, T]$, 

\begin{equation}
\label{x}
X(t)=\sigma \left(\int_{-\infty} ^ {t} e ^ {-\lambda (t-u) } du \right)(Z ^ {2}-1) = \frac{\sigma }{\lambda }(Z ^ {2}-1).
\end{equation}

\end{prop}
{\bf Proof: }  We  prove the convergence of the finite dimensional distributions of $ X ^ {H}$ to those of $X$ when $H\to 1$. Consider $\alpha _{1},.., \alpha _{d} \in \mathbb{R}$ and $t_{1},.., t_{d}\in [0,T]$.  We prove that
$$\sum_{j=1}^ {d} \alpha_{j}X ^ {H}(t_{j}) \xrightarrow[H\to 1]{(d)}   \sum_{j=1}^ {d} \alpha_{j}X (t_{j}).$$
Notice that, by linearity
$$\sum_{j=1}^ {d} \alpha_{j}X ^ {H}(t_{j}) =\int_{\mathbb{R}} g(u)dZ^ {H}(u) $$
with
\begin{equation}
\label{g2}
g(u)=\sigma \sum_{j=1}^ {d} \alpha_{j} e ^ {-\lambda (t_{j}-u) }1_{(-\infty, t_{j})} (u).
\end{equation}

We need to show (\ref{cH1}).  First, notice that $g$ belongs to $L ^{1} (\mathbb{R})$ because
\begin{equation}\label{8m-2}
\int_{\mathbb{R}} \vert g(u)\vert du \leq \sigma \sum_{j=1} ^{d} \vert \alpha _{j}\vert \int_{-\infty} ^{t_{j}} e ^ {-\lambda (t_{j}-u) }du =\sigma \sum_{j=1} ^{d} \vert \alpha _{j}\vert\int_{0} ^{\infty} e ^{-\lambda u}du =\frac{\sigma}{\lambda} \sum_{j=1} ^{d} \vert \alpha _{j}\vert <\infty. 
\end{equation}
Concerning the first part of (\ref{cH1}), we have, with $g$ from (\ref{g2}),
\begin{eqnarray*}
&&\int_{\mathbb{R}}\int_{\mathbb{R}} dudv \vert g(u) g(v)\vert \vert u-v\vert ^{2\eps -1} \leq \sigma ^{2}\sum_{j,k=1} ^{d} \vert \alpha _{j} \alpha _{k}\vert  \int_{-\infty} ^{t_{j}} du \int_{-\infty}^{t_{k}} dv e ^ {-\lambda (t_{j}-u) }e ^ {-\lambda (t_{k}-v) }\vert u-v\vert ^{2\eps -1}\\
&\leq &
\sigma ^{2}\sum_{j,k=1} ^{d} \vert \alpha _{j} \alpha _{k}\vert  
\int_{0} ^{\infty}du \int_{0} ^{\infty} dv e ^{-\lambda (u+v)} \vert u-v -(t_{j}-t_{k})\vert ^{2\eps -1}.
\end{eqnarray*}
The fact that the above integral is finite follows from the computation of the term $I$ given by (\ref{15d-4}) below.

  Again the tightness is obtained since $X^{H}$ obviously satisfies (\ref{7d-1}). \qed

\subsubsection{Convergence when $H\to \frac{1}{2}$ }

We will denote by $(X_{0} (t))_{t\geq 0}$ the stationary Ornstein-Uhlenbeck process. Recall that the stationary Ornstein-Uhlenbeck process is obtained  from (\ref{y0}) by taking the initial condition $\xi= \sigma \int_{-\infty}^ {0} e ^{-\lambda u} dW(u)$  where $(W(u))_{u\in \mathbb{R}}$ is a Wiener process on the whole real line. Thus
\begin{equation}
\label{xo}
X _{0} (t)= \sigma \int_{-\infty}^ {t} e ^ {-\lambda (t-u)}d W(u) 
\end{equation}
with $\sigma, \lambda >0.$  The process $(X_{0}(t))_{t\geq 0} $ is a  centered  Gaussian process, with stationary increments, with covariance function
\begin{equation}\label{12d-1}
\mathbf{E} X_{0}(t) X_{0} (s) = \frac{\sigma ^ {2}}{2\lambda } e ^ {-\lambda \vert t-s\vert } 
\end{equation}
for every  $s,t \geq 0$. It follows from (\ref{12d-1}) that $X_{0}$ is a stationary Gaussian   process with stationary increments.

\begin{prop}
The process $(X ^ {H}(t))_{t\in [0, T]}$ converges weakly, in the space of continuous functions $C[0,T]$, to the stationary Ornstein-Uhlenbeck process $(X_{0} (t))_{t\in [0, T]}$ .
\end{prop}
{\bf Proof: } Consider $\alpha_{1},.., \alpha_{d}\in \mathbb{R}$ and $t_{1},.., t_{d}\in [0,T]$.  In order to show that the random variable $\sum_{j=1} ^ {d}\alpha _{j} X ^ {H}(t_{j})$ converges in distribution to $\sum_{j=1} ^ {d} \alpha _{j} X _{0}(t_{j}) $ as $H\to \frac{1}{2}$,  we use Proposition \ref{p3}. In order to apply this result,  we first notice that the function $g$ given by (\ref{g2}) is in $\mathcal{H}_{H}$ (acutually  the computations below will show that it also belongs to $\vert \mathcal{H}_{H}\vert $).  We also need to verify (\ref{sf}) and (\ref{cH12}) in Proposition \ref{p3}. Let us start by checking (\ref{sf}). We have to prove that 
\begin{equation}
\label{15d-1}
\mathbf{E} \left( \sum_{j=1} ^ {d}\alpha _{j} X ^ {H}(t_{j})\right) ^{2} \xrightarrow[H \to \frac{1}{2}]{} \mathbf{E} \left( \sum_{j=1} ^ {d} \alpha _{j} X _{0}(t_{j}) \right) ^{2}.
\end{equation}
First, notice that by (\ref{12d-1}),
\begin{eqnarray*}
\mathbf{E} \left( \sum_{j=1} ^ {d} \alpha _{j} X _{0}(t_{j}) \right) ^{2}&=&\sum_{j,k=1} ^{d} \alpha_{j} \alpha _{k} \mathbf{E} X_{0}(t_{j})X_{0}(t_{k}) \\
&=& \frac{ \sigma ^{2} }{2\lambda } \sum_{j,k=1} ^{d} \alpha_{j} \alpha _{k} e ^{-\lambda \vert t_{j}-t_{k}\vert }.
\end{eqnarray*}
On the other hand, 
\begin{eqnarray}
\mathbf{E} \left( \sum_{j=1} ^ {d} \alpha _{j} X^{H}(t_{j}) \right) ^{2}&=&\sigma ^{2} H(2H-1) \sum_{j,k=1} ^{d} \alpha_{j} \alpha _{k}\int_{-\infty} ^{t_{j}} du \int_{-\infty}^{t_{k}} dv e ^{-\lambda (t_{j}-u)}e ^{-\lambda (t_{k}-v)} \vert u-v\vert  ^{2H-2}\nonumber \\ \nonumber
&=&\sigma ^{2} H(2H-1) \sum_{j,k=1} ^{d} \alpha_{j} \alpha _{k}\int_{0} ^{\infty} du \int_{0} ^{\infty}dv e ^{-\lambda u} e ^{-\lambda v} \vert u-v-(t_{j}-t_{k})\vert ^{2H-2}.\nonumber \\ 
\label{15d-3}
\end{eqnarray}
We need to compute the integral 
\begin{equation}
\label{15d-4}
I= H(2H-1)\int_{0} ^{\infty} du \int_{0} ^{\infty}dv e ^{-\lambda u} e ^{-\lambda v} \vert u-v-K\vert ^{2H-2}
\end{equation}
and we can assume by symmetry that $K\geq 0$.
We have
\begin{eqnarray*}
I&=& H(2H-1) \int_{0} ^{\infty} dv e ^{-\lambda v} \int_{v+K} ^{\infty } du e ^{-\lambda u} (u-v-K) ^{2H-2}\\
&&+ H(2H-1) \int_{0} ^{\infty} dv e ^{-\lambda v} \int_{0} ^{v+K} du e ^{-\lambda u} (v+K-u) ^{2H-2}\\
&:=& I_{1}+I_{2}.
\end{eqnarray*}
Let us regard the summand $I_{1}$.  By the change of variables $\tilde{u}= u-(v+K)$, 
\begin{eqnarray*}
I_{1}&=& H(2H-1) \int_{0} ^{\infty} dv e ^{-\lambda v}  \int_{0} ^{\infty} du e ^{-\lambda (u+v+k)} u^{2H-2} \\
&=&H(2H-1)e ^{-\lambda K}  \int_{0} ^{\infty} dv e ^{-2\lambda v}   \lambda ^{1-2H} \int_{0} ^{\infty} dz e ^{-z} z ^{2H-2} \\&=&H(2H-1)e ^{-\lambda K} \lambda ^{1-2H}\frac{1}{2\lambda} \int_{0} ^{\infty} dz e ^{-z} z ^{2H-2}=H(2H-1)e ^{-\lambda K} \frac{1}{2\lambda ^{2H}} \Gamma (2H)
\end{eqnarray*}
where $\Gamma $ is the gamma function. By integrating by parts 
$$H(2H-1) \int_{0} ^{\infty} dz e ^{-z} z ^{2H-2} \xrightarrow[H\to  \frac{1}{2}]{}1$$
so
\begin{equation}
\label{i1}
I_{1}\xrightarrow[H\to  \frac{1}{2}]{} \frac{e ^{-\lambda K}}{4\lambda}.
\end{equation}
Concerning the summand $I_{2}$, with $\tilde{u}= v+K-u$,
\begin{eqnarray*}
I_{2}&=& H(2H-1) \int_{0} ^{\infty} dv e ^{-\lambda v} \int_{0} ^{v+K} du e ^{-\lambda (v+K-u)} u ^{2H-2}\\
&=&H(2H-1) e ^{-\lambda K}\int_{0} ^{\infty}  du e ^{\lambda u} u ^{2H-2} \int_{ (u-K)\vee 0} ^{\infty} dv e ^{-2\lambda v}\\
&=& H(2H-1) e ^{-\lambda K} \int_{0} ^{K} du e ^{\lambda u} u ^{2H-2} \int_{0} ^{\infty} dv e ^{-2\lambda v}\\&&+  H(2H-1) e ^{-\lambda K} \int_{K}^{\infty} due ^{\lambda u} u ^{2H-2} \int_{u-K} ^{\infty} dv e ^{-2\lambda v}\\
&=&H(2H-1) e ^{-\lambda K} \frac{1}{2\lambda}\int_{0} ^{K} du e ^{\lambda u} u ^{2H-2} + H(2H-1) e ^{\lambda K} \frac{1}{2\lambda}\int_{K}^{\infty} due ^{-\lambda u} u ^{2H-2}.
\end{eqnarray*}
By using the integration by parts, we obtain
\begin{eqnarray*}
I_{2}&=&\frac{H}{2\lambda} e ^{-\lambda K} \left[ e ^{\lambda K} K ^{2H-1}-\lambda \int_{0} ^{K} e ^{\lambda u} u ^{2H-1}du\right]+ \frac{H}{2\lambda} e ^{\lambda K} \left[ -e ^{-\lambda K} K ^{2H-1}+\lambda \int_{K}^{\infty} e ^{-\lambda u} u ^{2H-1}du\right].
\end{eqnarray*}
Since 
\begin{equation*}
\lambda \int_{0} ^{K} e ^{\lambda u} u ^{2H-1}du \xrightarrow[H\to  \frac{1}{2}]{} e ^{\lambda K}-1\mbox{ and } \lambda \int_{K}^{\infty} e ^{-\lambda u} u ^{2H-1}du\xrightarrow[H\to  \frac{1}{2}]{} e ^{-\lambda K} 
\end{equation*}
we get
\begin{equation}
\label{i2}
I_{2} \xrightarrow[H\to  \frac{1}{2}]{} \frac{e ^{-\lambda K}}{4\lambda}.\end{equation}
From (\ref{i1}) and (\ref{i2}), the integral $I$ from (\ref{15d-4}) verifies
\begin{equation}
\label{ii}
I\xrightarrow[H\to  \frac{1}{2}]{} \frac{e ^{-\lambda K}}{2\lambda}.
\end{equation} 
Relation (\ref{ii}), together with (\ref{15d-3}), will imply that  (\ref{15d-1}), i.e. the assumption (\ref{sf}) is verified.

 Let us now check  the assumption (\ref{cH12}). With $g$ from (\ref{g2}),\begin{eqnarray*}
&&  \int_{\mathbb{R}}...\int_{\mathbb{R}} du_{1}...du_{4} \vert g (u_{1})...g(u_{m}) \vert \vert  u_{1}-u_{2} \vert  ^ {H-1 }.... \vert  u_{4}-u_{1} \vert  ^ {-H-1 }\\
&\leq  &\sum_{j_{1}, j_{2},.., j_{4}=1}^ {d} \vert \alpha _{j_{1}}...\alpha _{j_{4}}\vert \int_{-\infty} ^ {t_{j_{1}}} du_{1}....\int_{-\infty} ^ {t_{j_{4}}} du_{m} e ^ {-\lambda (t_{j_{1}}-u_{1})}....e ^ {-\lambda (t_{j_{4}}-u_{4})} \vert  u_{1}-u_{2} \vert  ^ {H-1 }.... \vert  u_{4}-u_{1} \vert  ^ {-H-1 }\\
&=&\sum_{j_{1}, j_{2},.., j_{4}=1}^ {d} \vert \alpha _{j_{1}}...\alpha _{j_{4}}\vert \int_{0} ^{\infty} du_{1}....\int_{0}^{\infty} du_{4} e ^ {-\lambda (u_{1}+..+u_{4})}\\
&&\times  \vert u_{1}-u_{2}-(t_{j_{1}}-t_{j_{2}})\vert ^ {H-1}... \vert u_{4}-u_{1}-(t_{j_{4}}-t_{j_{1}})\vert ^ {H-1}\\
&\leq &e ^{\frac{\lambda }{2}  ( \vert t_{j_{1}}-t_{j_{2}}\vert +...+ \vert t_{j_{4}}-t_{j_{1}}\vert ) }\sum_{j_{1}, j_{2},.., j_{4}=1}^ {d} \vert \alpha _{j_{1}}...\alpha _{j_{4}}\vert \int_{0} ^{\infty} du_{1}...\int_{0} ^{\infty}du_{4}\nonumber \\
&& e ^{-\frac{\lambda }{2}(\vert u_{1}-u_{2}-(t_{j_{1}}-t_{j_{2}})\vert+...+ \vert u_{4}-u_{1}-(t_{j_{4}}-t_{j_{1}})\vert)}\\
&&\times 
\left( 1\vee \vert u_{1}-u_{2}-(t_{j_{1}}-t_{j_{2}})\vert ^{H-1}\right) ....\left( 1\vee \vert u_{4}-u_{1}-(t_{j_{4}}-t_{j_{1}})\vert ^{H-1}\right)
\end{eqnarray*}
Now, we will consider the set $T'$ of affine functionals on $\mathbb{R} ^ {4}$ given by
$$T' =\{  u_{1}-u_{2}-(t_{j_{1}}-t_{j_{2}}),...,  u_{4}-u_{1}-(t_{j_{4}}-t_{j_{1}}) \}.$$
As before, $T'$ is the only paddet subset of $T'$. 

We apply the power counting theorem with 
$$(\alpha_{1},..,\alpha _{4})= ( H-1,.., H-1) \mbox{ and } ( \beta_{1},..,\beta_{4})=(-\gamma, ..., -\gamma)$$
with $\gamma \in (1-\frac{1}{4}, 1)=(\frac{3}{4}, 1)$. We have
$$d_{0}(T')= 4-1+ 4(H-1)>0 \mbox{ if } H>\frac{1}{4}$$

and 
$$d_{\infty} (\emptyset)=4-1-4\gamma <0 \mbox{ if }\gamma >1-\frac{1}{4}=\frac{3}{4}.$$
Therefore, the function 
$$H\to  \int_{\mathbb{R}}...\int_{\mathbb{R}} du_{1}...du_{4} \vert g (u_{1})...g(u_{m}) \vert \vert  u_{1}-u_{2} \vert  ^ {H-1 }.... \vert  u_{4}-u_{1} \vert  ^ {H-1 }$$
is finite and continuous on the set $D=\{ H \in (0,1],  H>\frac{1}{4}\}$  which implies that condition (\ref{cH12}) is satisfied. The conclusion follows from Proposition \ref{p3}.\qed 

\section{Appendix: Multiple stochastic integrals and the Fourth Moment Theorem}
Here, we shall only recall some elementary
facts; our main reference is  \cite{N}. Consider
${\mathcal{H}}$ a real separable infinite-dimensional Hilbert space
with its associated inner product ${\langle
.,.\rangle}_{\mathcal{H}}$, and $(B (\varphi),
\varphi\in{\mathcal{H}})$ an isonormal Gaussian process on a
probability space $(\Omega, {\mathfrak{F}}, \mathbb{P})$, which is a
centered Gaussian family of random variables such that
$\mathbf{E}\left( B(\varphi) B(\psi) \right) = {\langle\varphi,
\psi\rangle}_{{\mathcal{H}}}$, for every
$\varphi,\psi\in{\mathcal{H}}$. Denote by $I_{q}$ the $q$th multiple
stochastic integral with respect to $B$. This $I_{q}$ is actually an
isometry between the Hilbert space ${\mathcal{H}}^{\odot q}$
(symmetric tensor product) equipped with the scaled norm
$\frac{1}{\sqrt{q!}}\Vert\cdot\Vert_{{\mathcal{H}}^{\otimes q}}$ and
the Wiener chaos of order $q$, which is defined as the closed linear
span of the random variables $H_{q}(B(\varphi))$ where
$\varphi\in{\mathcal{H}},\;\Vert\varphi\Vert_{{\mathcal{H}}}=1$ and
$H_{q}$ is the Hermite polynomial of degree $q\geq 1$ defined
by:\begin{equation}\label{Hermite-poly}
H_{q}(x)=(-1)^{q} \exp \left( \frac{x^{2}}{2} \right) \frac{{\mathrm{d}}^{q}%
}{{\mathrm{d}x}^{q}}\left( \exp \left(
-\frac{x^{2}}{2}\right)\right),\;x\in \mathbb{R}.
\end{equation}

We  recall the hypercontractivity property of multiple stochastic integrals. If $Y=I_{2}(f)$, with $f\in \mathcal{H} ^{\otimes 2}$, then (see relation (2.7.2) in \cite{NPbook})
\begin{equation}
\label{hyper}
\mathbf{E} \left| Y\right| ^{q}\leq (q-1) ^{q} E \vert Y^2\vert ^{\frac{q}{2}}
\end{equation}
for every $q>2$.
We will use the following  famous result initially proven in \cite{NuPe} that characterizes the convergence in distribution of a sequence of multiple integrals torward the Gaussian law.

\begin{theorem}\label{fmt}
Fix $n \geq 2$ and let $\left( F_{k} , k \geq 1 \right)$ , $F_{k} = I_{n}\left(f_{k} \right) $ (with $f_{k} \in {\mathcal{H}}^{\odot n}$ for every $k \geq 1$), be a sequence of square-integrable random variables in the nth Wiener chaos such that $ \mathbf{E} \left[ F_{k}^{2} \right] \rightarrow 1$ as $k \rightarrow \infty$. The following are equivalent:
\begin{enumerate}
\item the sequence $\left( F_{k} \right)_{k \geq 0}$ converges in distribution to the normal law $\mathcal{N} (0,1)$; 
\item $ \mathbf{E} \left[ F_{k}^{4} \right] = 3$ as $k \rightarrow \infty$;
\item for all  $1 \leq  l \leq n-1$, it holds that $\lim\limits_{k \rightarrow \infty} \Vert  f_{k}  \otimes_{l}  f_{k} \Vert_{{\mathcal{H}}^{\otimes 2 (n-l)}} = 0 $;

\end{enumerate}
\end{theorem}
Other equivalent condition can be stated in term of the Malliavin derivatives of $F_{k}$, see \cite{NPbook}.

\end{document}